\documentclass[12pt,twoside]{amsart}
\usepackage{amssymb}
\usepackage{amscd}
\usepackage[abbrev,alphabetic]{amsrefs}
\usepackage{hyperref}

\usepackage{xypic}
\usepackage[all,cmtip]{xy}

\usepackage{caption}

\usepackage{url}
\RequirePackage{booktabs, multirow}
\RequirePackage{pgf}

\title[On bpf theorem for 3-folds in large characteristic]
{On the base point free theorem for klt threefolds in large characteristic} 
\author{Fabio Bernasconi} 

\address{Department of Mathematics, University of Utah, Salt Lake City, UT 84112, USA} 
\email{fabio@math.utah.edu}

\newcommand{\red}[0]{{\operatorname{red}}}

\newcommand{\Spec}[0]{{\operatorname{Spec}}}

\newcommand{\im}[0]{{\operatorname{im}}}

\newcommand{\Ex}[0]{{\operatorname{Ex}}}

\newcommand{\Diff}[0]{{\operatorname{Diff}}}

\newtheorem{thm}{Theorem}[section]
\newtheorem{lem}[thm]{Lemma}
\newtheorem{cor}[thm]{Corollary}
\newtheorem{proposition}[thm]{Proposition}

\newtheorem{claim}[thm]{Claim}

\theoremstyle{definition}

\newtheorem{definition}[thm]{Definition}

\newtheorem{remark}[thm]{Remark}

\makeatletter
  
  \@addtoreset{equation}{section}
  \makeatother

\newcommand{\MO}{\mathcal{O}}

\newcommand{\Q}{\mathbb{Q}}

\begin{document}

\begin{abstract}
In this article we present a refinement of the base point free theorem for threefolds in positive characteristic. If $L$ is a nef Cartier divisor of numerical dimension at least one on a projective Kawamata log terminal threefold $(X, \Delta)$ over a perfect field $k$ of characteristic $p \gg 0$ such that $L-(K_X + \Delta)$ is big and nef, then we show that the linear system $|mL|$ is base point free for all sufficiently large integers $m>0$. 
\end{abstract}

\maketitle

\tableofcontents

\section{Introduction}
The base point free theorem is one of the cornerstones of the Minimal Model Program (MMP for short) over a field of characteristic zero.
Due to the failure of the Kodaira vanishing theorem and its generalisations, it is not known whether it still holds for varieties over fields of positive characteristic.

However, in the case of threefolds the base point free theorem has been established with increasing generality in recent years. 
Let $(X,\Delta)$ be a projective klt threefold pair over a perfect field of characteristic $p>0$ and let $L$ be a nef Cartier divisor such that $L-(K_X+\Delta)$ is big and nef.
In the seminal article \cite{Kee99}, Keel proved that if $L$ is assumed moreover to be big, then it is endowed with a map (EWM) without any restriction on the characteristic.
After the development of the MMP for threefolds (\cite{HX15}), in \cite[Theorem 1.1]{Xu15}, \cite[Theorem 1.5]{Bir16} and \cite[Theorem 1.2]{BW17} the authors prove that the linear system $|mL|$ is base point free for sufficiently large and \emph{sufficiently divisible} $m>0$ if the characteristic $p$ is at least five whether or not $L$ is big. 

Let us recall that over a field of characteristic zero, the divisibility condition on $m$ is indeed superfluous (see \cite[Theorem 3.3]{KM98}). 
Thus one may wonder whether we can remove it also in characteristic $p>0$.
Unfortunately, this is not possible in low characteristic: in \cite[Theorem 1.2]{Tana}, Tanaka showed that the divisibility assumption is indeed necessary over fields of characteristic two and three when the numerical dimension of $L$ is one.

The aim of this article is to present a refinement of the base point free theorem for threefolds by proving that we can remove the divisibility assumption if the characteristic is sufficiently large.

\begin{thm}[Theorem \ref{t-bpf-final}]\label{t-intro1}
There exists a constant $p_0>5$ such that the following holds.
Let $k$ be a perfect field of characteristic $p> p_0$.
Let $(X, \Delta)$ be a  projective klt threefold log pair over $k$. 
Let $L$ be a nef Cartier divisor on $X$ such that
\begin{enumerate}
\item the numerical dimension $\nu(L)$ is at least one;
\item $nL-(K_X+\Delta)$ is a big and nef $\Q$-Cartier $\Q$-divisor for some $n>0$.
\end{enumerate}
Then there exists an integer $m_0>0$ such that the linear system $|mL|$ is base point free for all $m \geq m_0$.
\end{thm}

\begin{remark}
The author does not know whether Theorem \ref{t-intro1} can be extended to the case where $\nu(L)=0$.
This is related to understanding the torsion of the Picard group of varieties of Fano type in positive characteristic, which has recently attracted the attention of various authors (see \cite{Tana,BT19, FS19}). 
\end{remark}

Apart from the intrinsic interest of understanding the differences between characteristic zero and characteristic $p$ birational geometry,
Theorem \ref{t-intro1} is important for effective statements in positive characteristic.
Indeed, if the divisibility required on $m$ is arbitrarily large, there is no hope that the Effective base point free theorem of Koll\'ar (see \cite{Kol93}) could hold in positive characteristic.

Theorem \ref{t-intro1} is a consequence of the following descent result for Cartier divisors which are relatively numerically trivial on $(K_X+\Delta)$-negative contractions, which is the main technical result of this paper.

\begin{thm}[Theorem \ref{t-thm-bir}]\label{t-main-thm}
There exists a constant $p_0>5$ with the following property.
Let $k$ be a perfect field of characteristic $p> p_0$.
Let $\pi \colon X \rightarrow Z$ be a projective contraction morphism between quasi-projective normal varieties over $k$.
Suppose that there exists an effective $\Q$-divisor $\Delta \geq 0$ such that
\begin{enumerate}
\item $(X, \Delta)$ is a klt threefold log pair;
\item  $-(K_X+\Delta)$ is $\pi$-big and $\pi$-nef;
\item $\dim(Z) \geq 1$.
\end{enumerate} 
Let $L$ be a Cartier divisor on $X$ such that $L \equiv_{\pi} 0$.
Then there exists a Cartier divisor $M$ on $Z$ such that $L \sim \pi^*M$.
\end{thm}

\begin{remark}
The constant $p_0$ in Theorem \ref{t-intro1} and Theorem \ref{t-main-thm} comes from the Kawamata-Viehweg vanishing theorem for surfaces of del Pezzo type in large characteristic (see \cite[Theorem 1.2]{CTW17}).
\end{remark}

\begin{remark}
In a recent preprint \cite{Tanb}, H. Tanaka investigates freeness for Cartier divisors up to $p$-power in positive characteristic. 
Using the MMP for threefolds, he then shows a similar result of Theorem \ref{t-main-thm} for nef Cartier divisors on klt threefold log pairs over perfect fields of characteristic $p>5$ up to taking $p$-powers (see \cite[Theorem 1.7]{Tanb}).
\end{remark}

\subsection{Sketch of the proof}
The proof of Theorem \ref{t-main-thm} is divided into two steps: first we discuss Cartier divisors which are numerically trivial for pl-contractions over surfaces and threefolds (see Theorem \ref{p-lb-plneg}).
Then we prove the general case (see Theorem \ref{t-thm-bir}).

Let us overview the case of pl-contractions treated in Section \ref{s-num-trv-pl}.
Let $(X, \Delta)$ be a $\Q$-factorial dlt threefold pair and let $S$ be a prime divisor in $\lfloor{ \Delta \rfloor}$.
Let $\pi \colon (X, \Delta) \to Y$ be a $(K_X+\Delta)$-negative contraction where $S$ is a prime divisor which is $\pi$-anti-nef and $\dim(Y) \geq 2$. 
Let $L$ be a $\pi$-numerically trivial Cartier divisor. 
We aim to prove that $L$ is $\pi$-trivial over a neighbourhood of $\pi(S)$.
Since $L|_S$ is $\pi|_S$-trivial (Proposition \ref{p-dlt-bpf}), it is sufficient to lift sections from $S$ to prove that $L$ is $\pi$-trivial.
In order to do so, we show that the higher direct image $R^1 \pi_* \MO_X(L-S)$ vanishes.
First we prove the vanishing in the case where the fibres of $\pi$ are at most one dimensional (Proposition \ref{p-descent-small}), for which we generalise a result of Das and Hacon (see Proposition \ref{p-vanishing-small}).
To prove the general case of pl-contractions in Theorem \ref{p-lb-plneg}, we use some techniques developed by Hacon and Witaszek to prove the rationality of klt threefold singularities (see \cite{HW17}). 
The main ingredients are the Kawamata-Viehweg vanishing theorem for surfaces of del Pezzo type in large characteristic (see \cite{CTW17}) and the aforementioned case of pl-contractions whose maximum dimension of the fibres is one.

In Section \ref{sec-bir-gen}, we prove Theorem \ref{t-main-thm}.
First in Subsections \ref{ss-desc-bir}, \ref{ss-desc-conic} we discuss the case where $\dim(Z) \geq 2$. 
The idea is to use the MMP and by replacing a fibre of $\pi \colon X \to Z$ with a surface of del Pezzo type by Proposition \ref{p-plt-blowup}, we can apply Theorem \ref{p-lb-plneg} to conclude the descent.
For the case where $\dim(Z)=1$ treated in Subsection \ref{ss-final-sec}, we blend the previous cases with some results on del Pezzo fibrations proven in \cite{BT19}.
\medskip

\indent \textbf{Acknowledgements:} 
I would like to thank P. Cascini, I. Cheltsov, C.D. Hacon, Z. Patakfalvi, H. Tanaka and J. Witaszek for useful discussions and comments on this article.
I am also thankful to the referee for reading carefully the manuscript and comments. 
Part of this work was done while I was visiting EPFL in March 2019 and I would like to thank Z. Patakfalvi for his generous hospitality.
The author was supported by the Engineering and Physical Sciences Research Council [EP/L015234/1].
\section{Preliminaries}

\subsection{Notation}
\begin{enumerate}
\item We will freely use the notation and terminology in \cite{Har77} 
and \cite{Kol13}. 
\item 
We say that a scheme $X$ is {\em excellent} (resp. {\em regular}) 
if 
the local ring $\MO_{X, x}$ at any point $x \in X$ is excellent (resp. regular). 
For the definition of excellent local rings, 
we refer to \cite[\S 32]{Mat89}.
\item Given a scheme $X$, we denote by $X_{\text{red}}$ the reduced subscheme of $X$ such that the induced morphism $X_{\text{red}} \to X$ is
surjective. 
\item For an integral scheme $X$, 
we define the {\em function field} $K(X)$ of $X$ 
to be $\MO_{X, \xi}$ for the generic point $\xi$ of $X$. 
\item 
For a scheme $B$, 
we say that $X$ is a {\em variety over} $B$ or a $B$-{\em variety} if 
$X$ is an integral scheme that is separated and of finite type over $B$. 
We say that $X$ is a {\em curve} over $B$ or a $B$-{\em curve} 
(resp. a {\em surface} over $B$ or a $B$-{\em surface}, 
resp. a {\em threefold} over $B$) 
if $X$ is a $B$-variety of dimension one (resp. two, resp. three).
\item Let $B$ be a regular excellent scheme.  We say that $(X, \Delta)$ is a \emph{log pair} if $X$ is a normal quasi-projective $B$-variety, $\Delta$ is an effective $\Q$-divisor and $K_{X/B}+\Delta$ is $\Q$-Cartier. We will use the terminology of \cite{Kol13} for singularities of pairs. 
We say that $\Delta$ is a \emph{boundary} if $\lfloor{ \Delta \rfloor}$ is a reduced effective divisor.
\item We say that a normal scheme $X$ is \emph{$\Q$-factorial} if every Weil divisor $D$ is $\Q$-Cartier (i.e. there exists an integer $m$ such that $mD$ is a Cartier divisor).
\item For an $\mathbb{F}_p$-scheme $X$ we denote by $F \colon X \to X$ the {\em absolute Frobenius morphism}. For a positive integer $e$ we denote by $F^e \colon X \to X$ 
the $e$-th iterated absolute Frobenius morphism.
\item Let $\pi \colon X \to Y$ be a proper morphism between normal varieties over a field $k$. We say the morphism $\pi$ is a \emph{contraction} if $\pi_* \MO_X = \MO_Y$.
\item Given $\pi \colon X \to Y$ be a proper morphism between normal varieties over $k$, we denote by $\rho(X/Y)$ the \emph{relative Picard number}. 
If $\pi$ is a birational morphism, we denote by $\Ex(\pi)$ the exceptional locus of $\pi$.
\item A $k$-projective variety $X$ is said to be \emph{of Fano type} if there exists an effective boundary $\Delta$ such that $(X, \Delta)$ is klt and $-(K_X+\Delta)$ is ample. If $\dim(X)=2$ we say $X$ is \emph{of del Pezzo type}.
\item We refer to \cite{Laz04} for the definition of ampleness, nefness and bigness of Cartier divisors in the relative setting.
\end{enumerate}

\subsection{Numerically trivial Cartier divisors on excellent surfaces}

Even if we are interested in varieties over a perfect field $k$, we will frequently localise at non-closed points and we will need descent results for numerically trivial Cartier divisors on schemes essentially of finite type over a field $k$. 
For this reason we recall the base point free theorem for excellent surfaces proven by Tanaka (see \cite{Tan18a}).

\begin{proposition} \label{p-dlt-bpf}
Let $B$ be a regular excellent separated scheme of finite dimension.
Let $\pi \colon X \to S$ be a projective $B$-morphism of normal quasi-projective $B$-schemes.
Suppose $(X,\Delta)$ is a $\Q$-factorial $B$-surface log pair where $\Delta$ is a boundary.
Let $L$ be a $\pi$-nef Cartier divisor on $X$ and suppose that $L-(K_{X/B}+\Delta)$ is $\pi$-ample.
Then the following hold.
\begin{enumerate}
\item Suppose that $L \not \equiv_{\pi} 0$. Then there exists $b_0$ such that for all $b \geq b_0$, $b L$ is $\pi$-free. In particular, there exists a factorisation $\pi \colon X \xrightarrow{\pi_L} T \xrightarrow{g} S$ such that $L \sim \pi_L^* H$ for a $g$-ample Cartier divisor $H$ on $T$.
\item Suppose that $(X, \Delta)$ is dlt and $L \equiv_{\pi} 0$ and the Stein factorisation of $X \to S$ is $X \to \Spec(k) \to S$, where $k$ is a perfect field.
Then $L \sim 0$.
\end{enumerate}
\end{proposition}


\begin{proof}
Since $X$ is $\Q$-factorial, we can perturb the boundary $\Delta$ and assume that $\lfloor{ \Delta \rfloor}=0$.
We show (1).
The first part of the statement is \cite[Theorem 4.2]{Tan18a}.
As for the remaining part, there exist $H_{b_0}$ and $H_{b_0+1}$ Cartier divisors on $T$ such that $b_0 L = \pi_{L}^* H_{b_0}$ and $(b_0+1)L=\pi_{L}^* H_{b_0+1}$. 
In particular, $L= \pi_{L}^* (H_{b_0 +1} - H_{b_0}) $, thus concluding.

We show (2).
Since $\nu(L)=0$, then $(X, \Delta)$ is a log del Pezzo pair over a perfect field.
By the Riemann-Roch formula, we have $h^0(L) \geq 1$. Since $\nu(L)=0$, we conclude $L \sim 0$.
\end{proof}



\subsection{Pl-contractions}
We introduce the notion of (weak) pl-contractions, which is a natural generalisation of the notions of pl-divisorial contractions and pl-flipping contractions (see \cite[Definition 3.2]{GNT}) for the case of contractions of fibre type.

\begin{definition}\label{d-pl-cont}
Let $k$ be a field.
Let $(X, \Delta)$ be a dlt pair over $k$ and let $S$ be a prime divisor contained in $\lfloor{ \Delta \rfloor}$.
Let $\pi \colon X \to Y$ be a projective $k$-morphism between quasi-projective normal varieties.
We say that $\pi$ is a \emph{$(K_X+\Delta, S)$-pl-contraction} (resp. a \emph{weak $(K_X+\Delta, S)$-pl-contraction}) if
\begin{enumerate}
\item $-(K_X+\Delta)$ is $\pi$-ample,
\item $-S$ is $\pi$-ample (resp. $\pi$-nef).
\end{enumerate}
\end{definition}

We collect some properties of weak pl-contractions for later use.

\begin{lem} \label{l-fibers-pl-cont}
Let $k$ be a field.
Let $X$ be a normal variety over $k$ and let $S$ be a $\Q$-Cartier prime divisor.
Let $\pi \colon X \to Y$ be a proper contraction morphism between normal varieties over $k$ such that $-S$ is $\pi$-nef.
Then for all closed points $x \in \pi(S)$, we have $\pi^{-1}(x) \subset S$.
\end{lem}

\begin{proof}
Immediate since $-S$ is $\pi$-nef.
\end{proof}

\begin{lem}\label{l-semi-ampleness}
Let $k$ be a perfect field of characteristic $p>5$.
Let $(X, \Delta)$ be a $\Q$-factorial threefold dlt pair over $k$ and let $S$ be a prime divisor contained in $\lfloor{ \Delta \rfloor}$.
Let $\pi \colon X \to Y$ be a weak $(K_X+\Delta,S)$-pl-contraction.
Then $-S$ is $\pi$-semi-ample.
\end{lem}

\begin{proof}
We write 
$-S= K_X+\Delta - (K_X+\Delta+S)$.
Thus we conclude $-S$ is $\pi$-semi-ample by the relative base point free theorem (see \cite[Theorem 2.9]{GNT}).
\end{proof}

\begin{lem} \label{l-dim-fibers}
Let $k$ be a perfect field of characteristic $p>5$.
Let $(X, \Delta)$ be a $\Q$-factorial threefold dlt pair over $k$ and let $S$ be a prime divisor contained in $\lfloor{ \Delta \rfloor}$.
Let $\pi \colon X \to Y$ be a weak $(K_X+\Delta,S)$-pl-contraction.
Let $\pi \colon X \xrightarrow{g} Z \xrightarrow{h} Y$ be the relative semi-ample fibration associated to $-S$ given by Lemma \ref{l-semi-ampleness}.
Assume that
\begin{enumerate}
\item $-S$ is not $\pi$-ample over any neighbourhood of $\pi(S)$,
\item $\dim(Y) \geq 2$.
\end{enumerate} 
Then the dimension of the fibres of $g$ are at most one in a neighbourhood of $g(S)$.
\end{lem}

\begin{proof}
We can assume that $\pi(S)$ is a closed point $y \in Y$, as otherwise the conclusion is immediate. 
Since $\dim(Y) \geq 2$, there is a neighbourhood $U$ of $y$ such that the dimension of the fibres of $\pi$ over $U \setminus y$ is at most one.
If $S$ gets contracted to a point by $g$, we have that $h$ is an isomorphism over an open neighbourhood of $y$, 
thus contradicting assumption (1). 
\end{proof}

The following extraction result, which motivated the definition of weak pl-contraction, will be used repeatedly in the following sections.
It is essentially stating that given a $(K_X+\Delta)$-negative contraction from a threefold onto a variety of positive dimension we can replace, after some birational modification, a fibre with a surface of del Pezzo type.

\begin{proposition}[cf. {\cite[Proposition 2.15]{GNT}}]\label{p-plt-blowup}
Let $k$ be a perfect field of characteristic $p>5$.
Let $\pi \colon X \to Z$ be a projective contraction morphism of normal quasi-projective varieties over $k$ with the following properties:
\begin{enumerate}
\item [i)] $(X,\Delta)$ is a klt threefold log pair,
\item [ii)] $-(K_X+\Delta)$ is $\pi$-big and $\pi$-nef.
\item [iii)]  $0<\dim(Z) \leq 3$.
\end{enumerate}
Fix a closed point $z \in Z$. Then there exists a commutative diagram of quasi-projective normal varieties
\[
\begin{CD}
W @> \varphi >> Y \\
@V \psi VV @VV g V\\
X @> \pi >> Z,
\end{CD}
\]
and an effective $\Q$-divisor $\Delta_Y$ on $Y$ such that
\begin{enumerate}
\item $(Y, \Delta_Y)$ is a $\Q$-factorial plt pair;
\item $S=(g^{-1}(z))_{\red}$ is an irreducible component of $\lfloor{\Delta_Y \rfloor}$ and $g$ is a weak $(K_Y+\Delta_Y, S)$-pl-contraction;
\item $W$ is a smooth threefold and $\varphi$ and $\psi$ are projective birational morphisms.
\end{enumerate}
In particular, $S$ is a surface of del Pezzo type by adjunction.
\end{proposition}

\subsection{Vanishing for pl-contractions with one dimensional fibres}

In this subsection, we present a generalisation of a relative vanishing theorem of Kodaira type in positive characteristic due to Das and Hacon (see \cite{DH16}).
This will be used in Proposition \ref{p-descent-small} to prove a descent result for numerically trivial Cartier divisors on pl-contractions with one-dimensional fibres.
We refer to \cite{Sch14} for a thorough treatment of the trace map of the Frobenius morphism for log pairs.
\begin{proposition}[{cf. \cite[Proposition 3.2]{DH16}}]\label{p-F-spl-seq}
Let $k$ be a perfect field of characteristic $p>5$.
Let $(X, \Delta=S+\Gamma)$ be a $\Q$-factorial threefold plt pair over $k$ where $S$ is a prime Weil divisor and $\lfloor \Gamma \rfloor=0$.
Suppose that there exists an integer $f>0$ such that $(p^f-1)(K_X+\Delta)$ is an integral Weil divisor.
Then there exists $e>0$ such that the morphism induced by the trace map of Frobenius
\[ \psi_{ne} \colon F^{ne}_* \MO_X((1-p^{ne})(K_X+\Gamma) -p^{ne}mS) \to \MO_X(-mS).\]
is surjective for all $n \geq 1$ and for all $m \geq 1$ at all codimension two points of $X$ contained in $S$.
\end{proposition}
\begin{proof}
By the proof of \cite[Proposition 3.2]{DH16}, there exists an integer $e>0$ such that the natural morphism
\[ F^{ne}_* \MO_X((1-p^{ne})(K_X+\Gamma)) \to \MO_X \]
admits a splitting in the category of $\MO_X$-modules at all codimension two points of $X$ contained in $S$ for all $n \geq 1$ (cf. \cite[Equations 3.4, 3.5]{DH16}).
By tensoring with $\mathcal{O}_X(-mS)$ and considering reflexive hulls, we thus deduce that
\[ F^{ne}_* \MO_X((1-p^{ne})(K_X+\Gamma)-p^{ne}mS) \to \MO_X(-mS) \]
is surjective at all codimension two points of $X$ contained in $S$.
\end{proof}

The following vanishing theorem is an easy generalisation of \cite[Theorem 3.5]{DH16}.

\begin{proposition}\label{p-vanishing-small}
Let $k$ be a perfect field of characteristic $p>5$.
Let $(X, \Delta)$ be a $\Q$-factorial threefold dlt pair over $k$ and let $S$ be a prime divisor contained in $\lfloor{ \Delta \rfloor}$.
Let $\pi \colon X \to Y$ be a projective contraction morphism between normal quasi-projective varieties such that
\begin{enumerate}
\item the maximum dimension of the fibres of $\pi$ is one;
\item $\pi$ is a weak $(K_X+\Delta,S)$-pl-contraction.
\end{enumerate}
Let $L$ be a Cartier divisor on $X$ such that $L$ is $\pi$-nef.
Then for all $m > 0$ we have
\[ R^1\pi_* \MO_X(L-mS)=0 \]
in a neighbourhood of $\pi(S)$.
\end{proposition}

\begin{proof}
We follow the proof of \cite[Theorem 3.5]{DH16} and we show how to adapt their arguments in order to prove our statement.

Let us write $\Delta=S+\Delta'$. Since $X$ is $\Q$-factorial we can slightly perturb the boundary $\Delta'$ in order to find a boundary $\Gamma$ such that $(X, S+\Gamma)$ is plt and $(1-p^f)(K_X+S+\Gamma)$ is an integral $\Q$-Cartier Weil divisor for some $f>0$ and $\pi$-ample.
 
If $\pi$ is birational, then by Lemma \ref{l-fibers-pl-cont} we can suppose $\Ex(\pi) \subset S$.
If $\dim(Y)=2$, then by Lemma \ref{l-fibers-pl-cont} $S$ is a vertical divisor for $\pi$ (i.e. $\pi(S)$ is an irreducible curve in $Y$).
Moreover, by localising at codimension one points of $Y$ and applying the relative Kawamata-Viehweg vanishing theorem for excellent surfaces (see \cite[Theorem 3.3]{Tan18a}),
we can suppose that $R^1\pi_* \MO_X(L-mS)$ is supported on a finite number of closed points of $\pi(S)$.

Consider now $e>0$ given by Proposition \ref{p-F-spl-seq}.
Using the trace map of the Frobenius morphism for all $n \geq 1$ we have the exact sequence
\begin{small}
\[ 0 \to \ker(\psi_{ne}) \to F^{ne}_* \MO_X((1-p^{ne})(K_X+\Gamma)+ p^{ne}L-p^{ne}mS) \xrightarrow{\psi_{ne}} \MO_X(L-mS). \]
\end{small}
Now let us split the sequence into two short exact sequences:
\begin{small}
\begin{equation}\label{1st-ex-seq}
 0 \to \ker(\psi_{ne}) \to F^{ne}_* \MO_X((1-p^{ne})(K_X+\Gamma)+ p^{ne}L -p^{ne}mS) \to \im(\psi_{ne}) \to 0;
\end{equation}
\end{small}
and
\begin{equation}\label{2nd-ex-seq}
 0 \to   \im(\psi_{ne})  \to  \MO_X(L-mS) \to \mathcal{G}_{ne} \to 0.
\end{equation}
Consider the long exact sequence in cohomology obtained by applying the push-forward $\pi_*$ to the short exact sequence (\ref{1st-ex-seq}).
Since the fibres of $\pi$ are at most one dimensional, we have $R^2 \pi_*\ker(\psi_e)=0$.
Let us denote $H=-(K_X+mS-L+\Gamma)=-(K_X+S+\Gamma)+ (L-(m-1)S)$ and let $r$ be the Cartier index of $H$.
Since $L$ and $-S$ are $\pi$-nef, we have $H$ is $\pi$-ample for all $m \geq 0$.
Write $(p^{ne}-1)=r \cdot a+b$ for some positive integers $a$ and $b$ such that $0 \leq b <r$.
Thus we have
\begin{multline*}
R^1\pi_* (F^{ne}_* \MO_X((1-p^{ne})(K_X+\Gamma)+p^{ne}L-p^{ne}mS))= \\
=F^{ne}_* (R^1\pi_* \MO_X(raH-b(K_X+\Gamma+mS-L)+L-mS)),
\end{multline*}
which vanishes for sufficiently large $n>0$ by Serre vanishing.
Therefore we conclude $R^1\pi_*\im(\psi_{ne})=0$.

Since $\psi_{ne}$ is surjective at codimension two points contained in $S$ by Proposition \ref{p-F-spl-seq}, we have $\mathcal{G}_{ne}$ is supported on a finite number of points.
Therefore $R^1\pi_* \mathcal{G}_{ne}=0$ and thus we conclude $R^1\pi_*\MO_X(L-mS)
=0$.
\end{proof}
\section{Numerically trivial Cartier divisors on pl-contractions} \label{s-num-trv-pl}
In this section, we prove a descent result for numerically trivial Cartier divisors on threefolds under weak pl-contractions over surfaces and threefolds (see Theorem \ref{p-lb-plneg}).

First, in Subsection \ref{ss-small} we discuss the case of a weak pl-negative contraction with fibres whose maximum dimension is one.
In this case the main tool we use is the vanishing theorem proven in Proposition \ref{p-vanishing-small}. 

In Subsection \ref{ss-gen} we discuss the general case. 
Our strategy is based on techniques developed by Hacon and Witaszek in \cite{HW17} to prove that klt threefold singularities are rational in large characteristic.
The main ingredients in the proof of Theorem \ref{p-lb-plneg} are Proposition \ref{p-descent-small} and the Kawamata-Viehweg vanishing theorem for log del Pezzo surfaces in large characteristic (\cite[Theorem 1.2]{CTW17}).

\subsection{Descent for pl-contractions with one dimensional fibres} \label{ss-small}

We apply the vanishing theorem of Proposition \ref{p-vanishing-small} to discuss descent of numerically trivial Cartier divisors in the case of pl-contractions with one-dimensional fibres.

\begin{proposition} \label{p-descent-small}
Let $k$ be a perfect field of characteritic $p>5$.
Let $(X, \Delta)$ be a $\mathbb{Q}$-factorial dlt threefold pair over $k$ and let $S$ be a prime divisor contained in $\lfloor{ \Delta \rfloor}$.
Let $\pi \colon X \rightarrow Y$ be a projective morphism between quasi-projective normal varieties over $k$ such that
\begin{enumerate}
\item the maximum dimension of the fibres of $\pi$ is one;
\item $\pi$ is a weak $(K_X+\Delta,S)$-pl-contraction.
\end{enumerate}
Let $L$ be a Cartier divisor on $X$ such that $L \equiv_{\pi} 0$.
Then $L \sim_{\pi} 0$ over a neighbourhood of $\pi(S)$.
\end{proposition}

\begin{proof}
Since the question is local over the base, we can assume $Y$ to be affine. 
By \cite[Theorem 2.11]{GNT} $S$ is a normal variety.
Moreover by adjunction the pair $(S, \Diff_S(\Delta'))$ is dlt where 
\[ (K_X + \Delta)|_{S} = K_{S} + \Diff_S(\Delta'), \]
where $\Delta'=\Delta-S$.
Let us denote by $\pi|_S \colon S \to T$ the restricted morphism, where $T:= \pi(S)$.
Since $-(K_S+\Diff_S(\Delta'))$ is $\pi|_S$-ample, Proposition \ref{p-dlt-bpf} implies $L|_{S} \sim_{\pi|_S} 0$.

By Lemma \ref{l-fibers-pl-cont} for any $y \in T$ we have $\pi^{-1}(y) \subset S$.
Since $L|_{S} \sim_{\pi|_S} 0$, it is sufficient to prove that 
the morphism
\[ \pi_* \mathcal{O}_X(L) \to (\pi|_S)_* \mathcal{O}_S(L) \] 
is surjective to prove that $L \sim_{\pi} 0$ over a neighbourhood of $T$.
By Proposition \ref{p-vanishing-small} we have $R^1 \pi_* \MO_X(L-S)=0$ and thus we conclude.
\end{proof}

\subsection{Descent for pl-contractions over threefolds and surfaces} \label{ss-gen}

We begin by recalling the Kawamata-Viehweg vanishing theorem for log del Pezzo surfaces in large characteristic (see \cite[Theorem 1.2]{CTW17}), which plays a key role in the proof of Theorem \ref{p-lb-plneg}.

\begin{thm} \label{c-KVV-dlt}
There exists a constant $p_0>0$ with the following property.
Let $k$ be a perfect field of characteristic $p>p_0$.
Let $(S, \Delta)$ be a dlt surface log pair over $k$ such that $-(K_S+ \Delta)$ is ample.
Let $B$ be an effective $\Q$-divisor such that $(S, B)$ is dlt and let $L$ be a Weil divisor such that $L-(K_S+B)$ is ample.
Then $H^1(S, L) =0$.
\end{thm}

\begin{proof}
Fix $p_0$ as in \cite[Theorem 1.2]{CTW17}.
By a perturbing the coefficients of $\Delta$ and $B$, there exists $\Delta'$ (resp. $B'$) such that $(S, \Delta')$ (resp. $(S, B')$) is klt and $-(K_S+ \Delta')$ (resp. $L-(K_S+B')$) is a $\Q$-Cartier $\Q$-ample divisor.
We conclude by \cite[Theorem 1.2]{CTW17}.
\end{proof}

\begin{remark}
We do not know an explicit bound on $p_0$.
However the examples constructed in \cite[Theorem 4.2]{CT19} and \cite[Theorem 1.1]{Ber} show that $p_0>3$.
\end{remark}

We recall the restriction short exact sequence constructed in \cite{HW17}.

\begin{proposition} \label{p-ses}
Let $k$ be a perfect field of characteristic $p>5$.
Let $(X, \Delta)$ be a $\Q$-factorial dlt threefold log pair defined over $k$.
Let $S$ be an irreducible component of $\lfloor{\Delta \rfloor}$ and let $\Delta=\Delta' + S$.
Let $D$ be a Weil divisor on $X$. Then for all $n \in \mathbb{Z}_{\geq 0}$ there exists a short exact sequence
\[ 0 \to \mathcal{O}_X(-(n+1)S+D) \to \MO_X(-nS+D) \to \MO_S(G_n) \to 0, \]
where $G_n \sim_{\Q} -nS|_S+D|_S-\Delta_n$ and $\Delta_n \leq \Diff_S(\Delta')$.
\end{proposition}
\begin{proof}
Since $X$ is $\Q$-factorial, the pair $(X, S)$ is plt and we can apply \cite[Corollary 3.7]{HW17}.
Since $\Delta_n \leq \Diff_S(0) \leq \Diff_S(\Delta')$ we conclude.
\end{proof}

We are now ready to prove the main result of this section.

\begin{thm} \label{p-lb-plneg}
There exists an integer $p_0>5$ such that the following hold.
Let $k$ be a perfect field of characteritic $p>p_0$.
Let $(X, \Delta)$ be a $\mathbb{Q}$-factorial dlt threefold log pair over $k$, and let $S$ be a prime Weil divisor contained in $\lfloor{ \Delta \rfloor}$.
Let $\pi \colon X \rightarrow Y$ be a projective contraction morphism between quasi-projective normal varieties over $k$ such that
\begin{enumerate}
\item $\pi$ is a weak $(K_X+\Delta,S)$-pl-contraction;
\item $\dim(Y) \geq 2$.
\end{enumerate}
Let $L$ be a Cartier divisor on $X$ such that $L \equiv_{\pi} 0$.
Then $L \sim_{\pi} 0$ over a neighbourhood of $\pi(S)$.
\end{thm}
\begin{proof}
Consider $p_0>5$ for which the statement of Theorem \ref{c-KVV-dlt} holds.
As in the proof of Proposition \ref{p-descent-small}, by adjunction we may write
\[ (K_X+\Delta)|_S=K_S+\Diff_S(\Delta') \]
where $\Delta'=\Delta-S$ and we have $L|_{S} \sim_{\pi|_S} 0$ by Proposition \ref{p-dlt-bpf}.
Since by Lemma \ref{l-fibers-pl-cont} for any $y \in T:=\pi(S)$ we have $\pi^{-1}(y) \subset S$, it is sufficient to prove that 
the morphism
\[ \pi_* \mathcal{O}_X(L) \to (\pi|_S)_* \mathcal{O}_S(L) \] 
is surjective.
To prove surjectivity, we show the vanishing of $R^1\pi_*\MO_X(L-S)$.

Let $n \geq 0$ and let us consider the following exact sequence given by Proposition \ref{p-ses}: 
\[ 0 \to \mathcal{O}_X(-(n+1)S) \rightarrow \mathcal{O}_X (-nS) \to \mathcal{O}_S (G_n) \to 0, \]
where $G_n \sim_{\Q} -nS|_S - \Delta_n$ for some $\Delta_n \leq \Diff_S(\Delta')$.
We tensor with $L$ and we consider the following exact sequence of $\MO_Y$-modules obtained by applying the push-forward $\pi_*$:
\begin{small}
\begin{align*} 
R^1 \pi_*\mathcal{O}_X(L-(n+1)S) \to  R^1 \pi_*\mathcal{O}_X(L-nS) \to  R^1 (\pi|_S)_* \mathcal{O}_S(L+G_n).
\end{align*}
\end{small}
Thus in order to prove $R^1\pi_* \mathcal{O}_X(L-(n+1)S)=0$ for all $n \geq 0$
it is sufficient to prove the following two vanishing results in cohomology:
\begin{enumerate}
\item [(i)] $R^1(\pi|_S)_* \mathcal{O}_S(L + G_n) =0$ for every $n \geq 0$.
\item [(ii)] $R^1 \pi_* \mathcal{O}_X (L -mS)=0$ for $m \gg 0$ and sufficiently divisible.
\end{enumerate} 
To prove (i) let us note that
\begin{equation*}
\begin{split}
L|_{S} + G_n & \sim_{\Q}  (K_{S}+ \Diff_S(\Delta') -\Delta_n) -(K_S+\Diff_S(\Delta'))+ L|_S -nS|_{S} \\ 
& \sim_{\Q}  (K_{S}+ B) +A,
\end{split}
\end{equation*}
where $(S, B:=\Diff_S(\Delta') -\Delta_n)$ is klt and $A$ is $(\pi_S)_*$-ample.
If $\dim(T)\geq 1$, we conclude by the relative Kawamata-Viehweg vanishing (see \cite[Proposition 3.2]{Tan18a}).
If $\dim(T)=0$, we have
\[R^1 (\pi|_S)_*\mathcal{O}_S(L+G_n)=H^1(S, \mathcal{O}_S(L+G_n)). \]
Since $(S, \Diff_S(\Delta'))$ is a dlt pair such that $-(K_S+\Diff_S(\Delta'))$ is an ample $\Q$-Cartier $\Q$-divisor, we conclude $H^1(S, \mathcal{O}_S(L+G_n))=0$ by Theorem \ref{c-KVV-dlt}.

We now prove (ii).
If $-S$ is $\pi$-ample over a neighbourhood of $T$ we conclude by the relative Serre vanishing theorem.
Thus we can suppose $-S$ is not $\pi$-ample  over any neighbourhood of $T$.
By Lemma \ref{l-semi-ampleness} we can consider the semi-ample fibration over $Y$ associated to $-S$
\[ \pi \colon X \xrightarrow{g} Z \xrightarrow{h} Y, \]
and let us consider an integer $k>0$ such that $-kS=g^*H$ for an $h$-ample Cartier divisor on $Z$.

Since by Lemma \ref{l-dim-fibers} the fibres of $g$ are at most one dimensional, we can apply Proposition \ref{p-vanishing-small} to deduce that $R^ig_* \MO_X(L-mS)=0$ for $i>0$ and $m>0$.
By Proposition \ref{p-descent-small} we have $L=g^*M$ for some Cartier divisor $M$ which is $h$-trivial.
By the Grothendieck spectral sequence, we thus deduce 
$R^1\pi_*\MO_X(L-mkS)=R^1 h_* \MO_Z(M+mH)$.
Since $H$ is $h$-ample, we may apply relative Serre vanishing to $R^1 h_* \MO_Z(M+mH)$ to conclude that for all $m$ sufficiently large we have $R^1\pi_*\MO_X(L-mkS)=0$.
\end{proof}

\section{The base point free theorem in large characteristic} \label{sec-bir-gen}
The aim of this section is to prove Theorem \ref{t-general-desc}.
For this we discuss first descent of numerically trivial Cartier divisors in the birational case (see Subsection \ref{ss-desc-bir}) and in the case of conic bundles (see Subsection \ref{ss-desc-conic}). 
In these cases, the main techniques used are the MMP for threefolds and Proposition \ref{p-plt-blowup} to reduce to the case of pl-contractions proven in Theorem \ref{p-lb-plneg}.

In subsection \ref{ss-final-sec}, we combine these results together with results on del Pezzo fibrations from \cite{BT19} to prove Theorem \ref{t-general-desc}.

\subsection{Birational case}\label{ss-desc-bir}

In this subsection, we prove descent of relatively numerically trivial Cartier divisors under $(K_X+\Delta)$-negative birational contraction morphisms of threefolds in large characteristic.

We need the following easy lemma on birational maps which are isomorphisms in codimension one.

\begin{lem} \label{l-pull-back-flips}
Let $k$ be a field. 
Let $f \colon X \to Z$ be a proper contraction morphism of normal varieties over $k$.
Let us consider the following commutative diagram
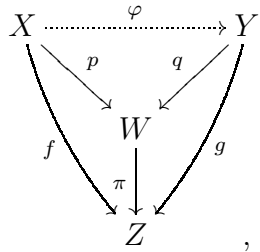
\[
\xymatrix{
&X \ar[dr]^p  \ar@/_/[ddr]_f \ar@{.>}[rr]^\varphi & & Y \ar[dl]_q \ar@/^/[ddl]^g \\
& & W \ar[d]_\pi & \\
& & Z & ,
}
\]
where $p \colon X \to W$ and $q \colon Y \to W$ are small proper birational contraction morphisms between normal varieties.
Let $L$ be a Cartier divisor on $X$ and suppose that there exists a Cartier divisor $M$ on $W$ such that $L \sim p^*M$.
Then $\varphi_*L $ is a Cartier divisor.
Moreover, the following are equivalent:
\begin{enumerate}
\item [(i)] there exists a Cartier divisor $H$ on $Z$ such that $L \sim f^*H$,
\item [(ii)]there exists a Cartier divisor $H$ on $Z$ such that $M \sim \pi^*H$,
\item [(iii)] there exists a Cartier divisor $H$ on $Z$ such that $\varphi_* L \sim g^*H$.
\end{enumerate}
\end{lem}
\begin{proof}
Since $\varphi_*L$ and $q^*M$ are both Weil divisors on a normal variety and they coincide outside a codimension two subset, we conclude that $\varphi_* L \sim q^*M$. 
In particular, $\varphi_* L $ is Cartier.

We now prove that (i) is equivalent to (ii). 
Obviously, (ii) implies (i).
If (i) holds, then $L \sim f^*H \sim p^*\pi^*H \sim p^*M$ implies $M \sim \pi^*H$. 
We can repeat the same proof using the equality $\varphi_*L \sim q^*M$ to conclude (ii) is equivalent to (iii).
\end{proof}

\begin{thm}\label{t-thm-bir}
There exists a constant $p_0>5$ with the following property.
Let $k$ be a perfect field of characteristic $p> p_0$.
Let $\pi \colon X \rightarrow Z$ be a projective birational morphism between quasi-projective normal varieties over $k$.
Suppose that there exists an effective $\Q$-divisor $\Delta$ such that
\begin{enumerate}
\item $(X, \Delta)$ is a klt threefold log pair,
\item $-(K_X+\Delta)$ is $\pi$-big and $\pi$-nef.
\end{enumerate} 
Let $L$ be a Cartier divisor such that $L \equiv_{\pi} 0$.
Then $L \sim_{\pi} 0$.
\end{thm}
\begin{proof}
Fix $p_0>5$ for which Theorem \ref{p-lb-plneg} holds.
Since the question is local over the base, we can  assume $Z$ to be affine and we fix $z \in Z$ a closed point.
By Proposition \ref{p-plt-blowup}, there exists a birational morphism
$g \colon Y \to Z$
such that there exists an effective $\Q$-divisor $\Delta_Y$ such that 
\begin{enumerate}
\item[(i)] $(Y, \Delta_Y)$ is a $\Q$-factorial plt pair,
\item[(ii)] $S:=(g^{-1}(z))_{\text{red}}$ is an irreducible component of $\lfloor{ \Delta_Y \rfloor}$ and $g$ is a weak $(K_X+\Delta, S)$-pl-contraction.
\end{enumerate}
Let us consider the following diagram
\[
\begin{CD}
W @> \varphi >> Y \\
@V \psi VV @VV g V\\
X @> \pi >> Z,
\end{CD}
\]
where $\varphi$ and $\psi$ are log resolutions.
Denote by $f:= \pi \circ \psi$.

In order to prove the theorem, it is sufficient to prove $\psi^*L \sim_f 0$.
We first prove that $\psi^*L$ descends to a Cartier divisor $M$ on $Y$.
To accomplish this, we apply Theorem \ref{p-lb-plneg} inductively.
\begin{claim}
There exists a Cartier divisor $M$ on $Y$ such that $\psi^*L \sim \varphi^*M$.
\begin{proof}
Since $Y$ is  $\Q$-factorial variety with klt singularities, we have
\[K_W + \sum_{i \in I} E_i \sim_{\Q} \varphi^*K_Y + \sum_{i \in I} (1+a_i) E_i,  \]
where $\bigcup_{i \in I} E_i = \text{Supp}(\text{Ex}(\varphi))$ and $1+a_i>0$ for all $i \in I$.
In particular,
\[ K_W +\sum_{i \in I} E_i \sim_{\Q, Y} \sum_{i \in I} (1+a_i) E_i.  \]
By \cite[Theorem 1.1]{HNT} we can run a $(K_W + \sum_{i \in I} E_i)$-MMP over $Y$
\[ h \colon W \dashrightarrow W_2 \dashrightarrow \cdots \dashrightarrow W_n=:T,\]
which terminates with a morphism $p \colon T \to Y$ such that $\sum_{i \in I}  (1+a_i) h_*E_{i}$ is $p$-nef. 
By the negativity lemma, this implies that all the divisors $E_i$ are contracted by $h$.
Thus $p \colon T \to Y$ is a small morphism. 
Since $Y$ is $\Q$-factorial, we conclude $p$ is an isomorphism.
Since we run a $(K_W+\sum_{i \in I} E_i)$-MMP over $Y$, we are also running a $\left(\sum_{i \in I} (1+a_i)E_i\right)$-MMP and thus each step is a pl-divisorial contraction or a pl-flip.
Thus we can apply Theorem \ref{p-lb-plneg} and Lemma \ref{l-pull-back-flips} inductively to conclude that there exists a Cartier divisor $M$ on $Y$ such that $\psi^*L \sim \varphi^*M$.
\end{proof}
\end{claim}
We can now apply Proposition \ref{p-lb-plneg} once more to show that there exists a Cartier divisor $N$ on $Z$ such that $M \sim g^*N$, thus concluding the proof.
\end{proof}

As a corollary, we obtain a descent result for numerically trivial Cartier divisors on threefolds admitting a birational morphism over a klt pair.
\begin{cor}\label{c-descent-base-klt}
There exists a constant $p_0>5$ with the following property.
Let $k$ be a perfect field of characteristic $p> p_0$.
Let $\pi \colon X \rightarrow Z$ be a projective birational morphism between quasi-projective normal varieties over $k$.
Suppose that there exists a $\Q$-divisor $B \geq 0$ such that $(Z,B)$ is a klt threefold.
Let $L$ be a Cartier divisor on $X$ such that $L \equiv_{\pi} 0$.
Then $L \sim_{\pi} 0$.
\end{cor}

\begin{proof}
Without any loss of generality, we can assume that $\pi \colon X \to Z$ is a log resolution for the pair $(Z, B)$.
Thus we can write 
\[ K_X + \pi_*^{-1}B+\sum_{i \in I}  E_i =\pi^*(K_Z+B) + \sum _{i \in I} (1+a_i) E_i. \]
Consider $0 < \varepsilon \ll 1$ such that $1+a_i-\varepsilon>0$ for any $i \in I$.
By \cite[Theorem 2.13]{GNT} we run a $(K_X + \pi_*^{-1}B+\sum_{i \in I} (1-\varepsilon) E_i)$-MMP over $Z$:
\[ f \colon X \dashrightarrow X_1 \dashrightarrow \dots \dashrightarrow X_n=:Y, \]
which ends with a relative minimal model $g \colon Y \to Z$.
By applying Lemma \ref{l-pull-back-flips} and Theorem \ref{t-thm-bir} at each step of the MMP inductively,
it is sufficient to prove that any Cartier divisor $N$ on $Y$ which is numerically $g$-trivial then it is $g$-trivial.
By the negativity lemma $g$ is a small birational morphism and thus we have $(Y, g_*^{-1} B)$ is klt and $(K_Y+g_*^{-1}B)$ is $g$-trivial.
In particular, $-(K_Y+g_*^{-1}B)$ is $g$-big and $g$-nef and thus we can apply Theorem \ref{t-thm-bir} once more to conclude.
\end{proof}

\subsection{Conic bundles} \label{ss-desc-conic}

In this section we prove descent of numerically trivial Cartier divisors under $(K_X+\Delta)$-negative contraction morphisms of relative dimension one (also known as conic bundles).

\begin{thm}\label{t-conic-descent}
There exists a constant $p_0>5$ with the following property.
Let $k$ be a perfect field of characteristic $p> p_0$.
Let $\pi \colon X \rightarrow Z$ be a projective contraction morphism between quasi-projective normal varieties over $k$.
Suppose that there exists an effective $\Q$-divisor $\Delta$ such that
\begin{enumerate}
\item $(X, \Delta)$ is a klt threefold log pair,
\item $-(K_X+\Delta)$ is $\pi$-big and $\pi$-nef,
\item $\dim(Z)=2$.
\end{enumerate} 
Let $L$ be a Cartier divisor on $X$ such that $L \equiv_{\pi} 0$.
Then $L \sim_{\pi} 0$.
\end{thm}

\begin{proof}
By Proposition \ref{p-dlt-bpf}, there exists open subset $U \subset Z$ such that $L|_{\pi^{-1}(U)} \sim_{\pi} 0$ and $Z \setminus U$ is a finite set of points.
Let $z$ be a closed point in $Z \setminus U$.
By Proposition \ref{p-plt-blowup},
there exists a birational morphism
$g \colon Y \to Z$
and an effective $\Q$-divisor $\Delta_Y$ on $Y$ such that 
\begin{enumerate}
\item[(i)] $(Y, \Delta_Y)$ is a $\Q$-factorial plt pair (in particular, $Y$ is klt),
\item[(ii)] $S:=(g^{-1}(z))_{\text{red}}$ is an irreducible component of $\lfloor{ \Delta_Y \rfloor}$ and $g$ is a weak $(K_X+\Delta, S)$-pl-contraction.
\end{enumerate}
Let us consider the following diagram
\[
\begin{CD}
W @> \varphi >> Y \\
@V \psi VV @VV g V\\
X @> \pi >> Z,
\end{CD}
\]
where $\varphi$ and $\psi$ are log resolutions.
Since $L$ is $\pi$-trivial in a punctured neighbourhood of $z$, 
to prove the statement it is sufficient to prove $\psi^*L$ is $(\pi \circ \psi)$-trivial over a neighbourhood of $z$. 
By Corollary \ref{c-descent-base-klt} there exists a Cartier divisor $M$ on $Y$ such that $\psi^*L = \varphi^*M$.
It is thus sufficient to prove that the Cartier divisor $M$ on $Y$ is $g$-trivial over a neighbourhood of $z$.
This is a consequence of Theorem \ref{p-lb-plneg}.
\end{proof}

\subsection{General case} \label{ss-final-sec}

We prove now the main theorem of this article. 
To deal with the remaining case of $(K_X+\Delta)$-negative contraction morphisms of relative dimension two, we combine Theorem \ref{t-thm-bir} and Theorem \ref{t-conic-descent} with a result on relatively numerically trivial Cartier divisors on del Pezzo fibrations (see \cite[Theorem 1.1]{BT19}), on the image of surfaces of del Pezzo type over imperfect fields of characteristic at least seven (see \cite[Corollary 5.8]{BT19}) and the MMP for $\Q$-factorial surfaces (see \cite{Tan18a}).

\begin{thm}\label{t-general-desc}
There exists a constant $p_0>5$ with the following property.
Let $k$ be a perfect field of characteristic $p>p_0$.
Let $(X, \Delta)$ be a klt threefold log pair and let $\pi \colon X \to Z$ be a projective contraction morphism between quasi-projective normal varieties over $k$.
Suppose that 
\begin{enumerate}
\item  $-(K_X+\Delta)$ is $\pi$-big and $\pi$-nef,
\item $\dim(Z) \geq 1$.
\end{enumerate} 
Let $L$ be a Cartier divisor such that $L \equiv_{\pi} 0$.
Then $L \sim_{\pi} 0$.
\end{thm}

\begin{proof}
We can assume $k$ is algebraically closed by a base change to an algebraic closure.
By taking a $\Q$-factoralization (\cite[Theorem 1.6]{Bir16}) we can further assume $X$ is $\Q$-factorial.
If $\dim(Z) \geq 2$, we conclude by Theorem \ref{t-thm-bir} and Theorem \ref{t-conic-descent}.

If $\dim(Z)=1$, by \cite[Theorem 2.12]{GNT} we can run a $(K_X+\Delta)$-MMP over $Z$:
\[ X_0:=X \dashrightarrow  X_1 \dashrightarrow \dots \dashrightarrow X_n =: Y, \]
which terminates with a Mori fibre space $g \colon Y \to T$ over $Z$ and let us denote by $f \colon Y \to Z$ the natural morphism. 
Let us note that since $X_{k(Z)}$ is a surface of del Pezzo type over $k(Z)$, we deduce that also the generic fibre $Y_{k(Z)}$ is a surface of del Pezzo type over $k(Z)$ by \cite[Lemma 2.9]{BT19}. 

By Corollary \ref{c-descent-base-klt} it is now sufficient to prove that given $M$ a Cartier divisor on $Y$ such that $M \equiv_{f} 0$, then $M \sim_{f} 0$.
We subdivide the proof according to the dimension of $T$.
If $\dim(T)=1$, then $T=Z$ and thus we conclude by \cite[Theorem 1.1]{BT19}.

If $\dim(T)=2$, we have a factorisation $Y \xrightarrow{g} T \xrightarrow{h} Z$ and by Theorem \ref{t-conic-descent} there exists a Cartier divisor $N$ on $T$ such that $N \equiv_{h}0$.
To conclude it is sufficient to prove that $N \sim_h 0$. 
Note that $T$ is a $\Q$-factorial surface by \cite[Theorem 5.4]{HNT}.
Since $Y_{k(Z)}$ is a surface of del Pezzo type, we deduce that the generic fibre $T_{k(Z)}$ is a Fano curve over $k(Z)$ by \cite[Corollary 5.8]{BT19}.
In particular $K_T$ is not pseudoeffective over $Z$.
Since $T$ is $\Q$-factorial, we can run a $K_T$-MMP over $Z$ by \cite[Theorem 1.1]{Tan18a}:
\[ \psi \colon T_0:=T \to T_1 \to T_2 \to \dots \to T_n=:V  \] 
which terminates with a Mori fibre space $p \colon V \to Z$.
By Proposition \ref{p-dlt-bpf} we show that there exists a Cartier divisor $D$ on $V$ such that $\psi^*D = N$.
Again by Proposition \ref{p-dlt-bpf} there exists a Cartier divisor $E$ on $Z$ such that $D=p^*E$, thus concluding that $N= h^*E$.
\end{proof}

We now show our improvement of the base point free theorem in large characteristic.

\begin{thm}\label{t-bpf-final}
There exists a constant $p_0>5$ such that the following holds.
Let $k$ be a perfect field of characteristic $p> p_0$.
Let $(X, \Delta)$ be a  quasi-projective klt threefold log pair and let $\pi \colon X \to Z$ be a projective contraction morphism of quasi-projective normal varieties over $k$.
Let $L$ be a $\pi$-nef Cartier divisor on $X$ such that
\begin{enumerate}
\item $\dim(Z) \geq 1$ or $\dim(Z)=0$ and $\nu(L) \geq 1$;
\item $nL-(K_X+\Delta)$ is a $\pi$-big and $\pi$-nef $\Q$-Cartier $\Q$-divisor for some $n>0$.
\end{enumerate}
Then there exists $m_0>0$ such that $mL$ is $\pi$-free for all $m \geq m_0$.
\end{thm}

\begin{proof}
By the relative base point free theorem for threefolds (see \cite[Theorem 2.9]{GNT}), $L$ is $\pi$-semi-ample.
Let $\pi \colon X \xrightarrow{f} Y \xrightarrow{g} Z $ denote the semi-ample fibration over $Z$ given by $L$. By hypothesis, we know that $\dim(Y) \geq 1$.
Since $L \equiv_{f} 0$, we deduce by Theorem \ref{t-thm-bir} that there exists a $g$-ample Cartier divisor $H$ on $Y$ such that $L \sim f^*H$.
This concludes the proof.
\end{proof}


\end{document}